\documentstyle[12pt]{article}

\newlength{\abstractwidth}
\renewcommand{\thefootnote}{\fnsymbol{footnote}}
\renewcommand{\thanks}[1]{\footnote{#1}} % Use this for footnotes
\newcommand{\starttext}{ \setcounter{footnote}{0}
\renewcommand{\thefootnote}{\arabic{footnote}}}

\newcommand{\be}{\begin{equation}}
\newcommand{\bea}{\begin{eqnarray}}
\newcommand{\eea}{\end{eqnarray}} \newcommand{\ee}{\end{equation}}

\def\ba{\begin{eqnarray}}
\def\ea{\end{eqnarray}}

\def\o{\omega}

\def\log{\,{\rm log}\,}
\def\exp{\,{\rm exp}\,}

\def\o{\omega}

\def\al{\alpha}

\def\o{\omega}

\def\vp{\varphi}

\def\na{\nabla}

\def\R{{\bf R}}

\def\i{\infty}

\def\ddb{{\partial\bar\partial}}

\def\na{{\nabla}}

\def\[{{\bf [}}
\def\]{{\bf ]}}

\def\pl{\partial}

\begin{document}
\starttext \baselineskip=15pt \setcounter{footnote}{0}
\newtheorem{theorem}{Theorem}
\newtheorem{lemma}{Lemma}
\newtheorem{definition}{Definition}
\newtheorem{proposition}{Proposition}

\begin{center}
{\Large \bf ON POINTWISE GRADIENT ESTIMATES FOR THE COMPLEX MONGE-AMPERE EQUATION
\footnote{Work supported in part by DMS-07-57372 and DMS-09-05873.}}
\bigskip\bigskip

{\large  D.H. Phong and Jacob Sturm} \\

\end{center}

\medskip

\begin{abstract}

In this note, a gradient estimate for the complex Monge-Amp\`ere equation is established.
It differs from previous estimates of Yau, Hanani, Blocki, P. Guan,
B. Guan - Q. Li in that it is pointwise, and  depends only on the
infimum of the solution instead of its $C^0$ norm.

\end{abstract}

\section{Introduction}
\setcounter{equation}{0}

A priori estimates for complex Monge-Amp\`ere equations are of considerable
interest in non-linear partial differential equations and complex geometry. Among
them, gradient estimates are somewhat special. In Yau's classic work on the
Calabi conjecture \cite{Y78}, they could be bypassed, since he was able to derive
a priori estimates for the Laplacian $\Delta\vp$ of the solution $\vp$, assuming only a
priori estimates for $\|\vp\|_{C^0}$. General linear elliptic theory allows then
to control $\|\na\vp\|_{C^0}$ in terms of $\|\Delta\vp\|_{C^0}$ and $\|\vp\|_{C^0}$.
Subsequently, several more direct estimates for $\|\na\vp\|_{C^0}$
were obtained by Hanani \cite{Hanani}, Blocki \cite{B08}, P. Guan \cite{G},
and B. Guan-Q. Li \cite{GL}
which did not require estimating $\|\Delta\vp\|_{C^0}$
first. This results in an improved dependence on the ambiant geometry, together with
greater flexibility for various generalizations, as we shall describe in more detail
in \S 2 below. In the case of the Dirichlet
problem, a different argument
was given by Chen \cite{C00} using blow-ups.

\medskip

The purpose of this note is to present a new gradient estimate, which has at least two
distinct advantages. The first is that it depends only on ${\rm inf}\,\vp$, and not
on $\|\vp\|_{C^0}$
as in all earlier gradient estimates. The second is that it is a more precise
pointwise estimate for $\nabla\vp(z)$ which remains valid even
when $\|\na\vp\|_{C^0}$ is unbounded. Such features were essential
in the construction in \cite{PS09} of geodesic rays starting from a test configuration,
and in fact, a version of the present gradient estimates had been established there.
That formulation was however somewhat obscured by the particular geometric set-up
of \cite{PS09}, and it seems worthwhile to extract a simpler and more general gradient
estimate, even though the proof is essentially the same, in the expectation that it will find other uses.

\section{The gradient estimate}

Let $(M,\omega)$ be a compact K\"ahler manifold with smooth boundary $\partial M$ (which may be empty)
and complex dimension $n$.
We consider the Monge-Amp\`ere equation on $\bar M$
\bea
\label{MA}
(\o+{i\over 2}\ddb \vp)^n=F(z,\vp)\,\o^n.
\eea
Here $F(z,t)$ is a  $C^2$ function on $\bar M\times\R$
which is assumed to be positive on the set $\bar M\times [\inf\,\vp,\infty)$.
We impose also the Dirichlet condition
\bea
\label{Dirichlet}
\vp=\vp_b\ \ {\rm on}\ \ \pl M\ \ {\rm if} \ \ \pl M\not=\emptyset
\eea
for a given $C^2$ function $\vp_b$.
Our main goal is to obtain gradient estimates with constants
which
depend only on the lower bound for $\vp$ (and not on ${\rm sup}_M|\vp|$).
There are two versions of such gradient estimates. In the first version,
the function $\vp$ is assumed to be $C^4$ on $\bar  M$. In the second,
$\vp$ is assumed to be $C^4$ on $ \bar  M\setminus S$, where $S$ is a set which
does not intersect $\pl M$, and $\vp(z)$ is assumed to tend to $\infty$
as $z\to S$. All covariant derivatives and curvatures listed below
are with respect to the metric $\o\equiv {i\over 2}g_{\bar kj}dz^j\wedge d\bar z^k$.

\begin{theorem}
\label{main}
Let $(M,\omega)$ be a compact K\"ahler manifold
with boundary $\pl M$ as described above, and let $F\in C^2(\bar M\times \R)$
satisfy $F>0$.

\smallskip

{\rm (a)} Let $\vp\in C^4(\bar M)$ be a solution
of the equation (\ref{MA}) (with the Dirichlet
condition (\ref{Dirichlet}) if $\partial M$ is not empty).
Then we have the a priori estimate
\bea
\label{C1a}
|\nabla \vp(z)|^2 \leq C_1\,{\rm exp}(A_1\, \vp(z)),
\qquad z\in {\bar M}
\eea
where $C_1$ and $A_1$ are constants that depend only on upper bounds on
\bea
\label{Lambda}
\Lambda=
-\inf_{a,b\in TM} {R_{\bar j i\bar l k}a^i\overline{a^j}b^k\overline{b^l}\over |a|^2|b|^2}
\eea
and 
\bea
\label{constants}
&&
{\rm inf}_M \vp,
\
{\rm sup}_{M\times [{\rm inf}\,\vp,\infty)} F,
\
{\rm sup}_{M\times [{\rm inf}\,\vp,\infty)}\left(|\nabla F^{1\over n}|+|\pl_tF^{1\over n}|\right),
\
\nonumber\\
&&
{\rm sup}_{\pl M}|\vp|,
\
{\rm sup}_{\pl M}|\nabla \vp|.
\eea

{\rm (b)} Let now $\vp\in C^4(\bar M\setminus S)$
be a solution of the equation (\ref{MA}) in $\bar M\setminus S$
with the Dirichlet condition (\ref{Dirichlet}). Assume further
that there exists a constant $B>0$ so that
\bea
\label{divisor}
&&
\vp(z)\to +\infty\ {\rm as}\ z\to\ S,\nonumber\\
&&
\log |\nabla \vp(z)|^2-B\,\vp(z)
\to -\infty \ {\rm as}\ z\to\ S.
\eea
Then we have
\bea
\label{C1b}
|\nabla \vp(z)|^2\leq C_1\, {\rm exp}(A_1\,\vp(z)),
\qquad z\in \bar M\setminus S,
\eea
where the constant $A_1={\rm max} (B,D)$, and the constants $C_1,D$
again depend only on the quantities listed in (\ref{Lambda}) and (\ref{constants}).

\end{theorem}

\bigskip

Before giving the proof of Theorem \ref{main}
we make a few remarks:

\smallskip
(1) The gradient estimates of Theorem \ref{main}
are the natural analogues of the classic estimates of Yau
and Aubin \cite{Y78,A} for the Laplacian $\Delta\vp$
of $\vp$. More precisely, under the same hypotheses as (a) of
Theorem \ref{main}, the estimates of Yau and Aubin are
\bea
\label{C2equation}
|\Delta\vp(z)|\leq C_2\exp({A_2(\vp(z)-\inf_M\vp)}),\quad z\in \bar M
\eea
where $A_2$ and $C_2$ are constants depending only on
upper bounds for $-\Delta\log F$,
the scalar curvature of $\o$, $\Lambda$,
$\|\vp\|_{C^0(\partial M)}$, and ${\rm sup}_{\partial M}(n+\Delta\vp)$.
On the other hand, when there is no boundary, and the function $F(z,\vp)$ is a function
$F(z)$ of $z$ alone, the equation
(\ref{MA}) is unchanged under shifts of $\vp$ by an additive constant.
Thus the infimum of $\vp(z)$ can be normalized to be $0$ by replacing $\vp(z)\to\vp(z)-{\rm inf}\vp$, so
we obtain the estimate
\bea
\label{C1elementary}
|\nabla\vp(z)|^2\leq C_1\,{\rm exp}(A_1(\vp(z)-\inf_M\,\vp)) \quad z\in \bar M
\eea
where the constant $C_1$ does not depend on ${\rm inf}\,\vp$, but depends only on the other quantities
in (\ref{constants}). Thus we see that the apriori estimate (\ref{C2equation}) for the Laplacian 
and the apriori estimate (\ref{C1elementary}) for the gradient have the same structure.

\medskip
(2) Not surprisingly, the constants
${\rm sup}_{M\times [{\rm inf}\,\vp,\i]} F$ and
${\rm sup}_{M\times [{\rm inf}\,\vp,\i]}|\nabla F^{1\over n}|+|\pl_t F^{1\over n}|$
in (\ref{constants}) can be replaced by
${\rm sup}_{M\times [{\rm inf}\,\vp,{\rm sup}\,\vp]} F$
and
${\rm sup}_{M\times [{\rm inf}\,\vp,{\rm sup}\vp]}|\nabla F^{1\over n}|+|\pl_t F^{1\over n}|$
respectively. Thus, when $\|\vp\|_{C^0}$ is bounded,
we obtain gradient bounds for $\vp$ for completely general
smooth and strictly positive functions $F(z,\vp)$.
We have however stated them in the original form since we are mainly
interested in the cases when there is no upper bound for ${\rm sup}\,\vp$.

\medskip

(3) As we had stressed,
the point of the above gradient estimates is that they depend
only on ${\inf}_M\vp$. If a dependence on $\|\vp\|_{C^0}$
is allowed, then there are many earlier direct approaches.
The first appears to be due to
Hanani \cite{Hanani}. More recently, Blocki \cite{B08}
gave a different proof, and our approach  builds directly on his.
The method of P. Guan \cite{G} can be extended to Hessian equations,
while the method of B. Guan-Q. Li \cite{GL} allows a general Hermitian metric $\omega$
as well as a more general right hand side $F(z)\chi^n$,
where $\chi$ is a K\"ahler form.

\bigskip

\noindent
{\it Proof of Theorem \ref{main}}: We adapt the proof from the earlier paper
\cite{PS09}. The key ingredient is an inequality due to Blocki \cite{B08}.

\smallskip
Let $g_{\bar kj}'=g_{\bar kj}+\pl_j\pl_{\bar k}\vp$.
We denote covariant derivatives with respect to $g_{\bar kj}'$ by $\na'$.
It is also convenient to formulate inequalities in terms of the relative
endomorphism
\bea
h^j{}_k=g^{j\bar p}g_{\bar p k}'.
\eea
Let $\gamma(x)$ be a function with $\gamma'(x)>0$ and $-\gamma''(x)>0$
for $x\in [{\rm inf}_M \vp,\infty)$. Set
\bea
\label{alpha}
\beta(z)=|\nabla \vp(z)|^2,
\qquad
\alpha(z)=\log\,|\nabla \vp(z)|^2-\gamma(\vp).
\eea
Then  Blocki shows that at an interior maximum for $\al$, the following
inequality holds:
\bea
\label{Blocki}
\Delta'\alpha
&\geq& {1\over \beta}|\nabla'\bar\nabla \vp|^2+
(\gamma'(\vp)-\Lambda -{F_1\over \beta^{1\over 2}}){\rm Tr}\,h^{-1}
+(- \gamma''(\vp)+2{\gamma'(\vp)\over \beta})|\nabla'\vp|^2
\nonumber\\
&&
-(n+2)\gamma'(\vp)-{2\over \beta}.
\eea
Here $\Lambda$ is the lower bound for the bisectional curvature defined in (\ref{Lambda}),
\bea
F_1=2\,{\rm sup}\,|\nabla(F(z,\vp(z))^{1\over n})|,
\eea 
and we have followed
\cite{PS09} in reformulating the inequality in terms of covariant derivatives
and the endomorphism $h$ and its inverse $h^{-1}$.

Since $\pl_z (F(z,\vp(z))^{1\over n})=(\pl_zF^{1\over n})(z,\vp)+\pl_z\vp(z)(\pl_t F^{1\over n})(z,\vp)$,
(\ref{Blocki}) implies
\bea
\label{Blocki1}
\Delta'\alpha
&\geq& {1\over \beta}|\nabla'\bar\nabla \vp|^2+
(\gamma'(\vp)-\Lambda -F_1'-{F_1''\over \beta^{1\over 2}}){\rm Tr}\,h^{-1}
+(- \gamma''(\vp)+2{\gamma'(\vp)\over \beta})|\nabla'\vp|^2
\nonumber\\
&&
-(n+2)\gamma'(\vp)-{2\over \beta}.
\eea
where we have set
\bea
F_1'={\rm sup}_{M\times [{\rm inf}\phi,\infty)}|\pl_t F^{1\over n}|,
\quad
F_1''={\rm sup}_{M\times [{\rm inf}\phi,\infty)}|\na_z F^{1\over n}|
\eea

\smallskip

We now prove Part (b) of Theorem \ref{main} (Part (a)
is an immediate consequence of Part (b), simply by taking $S$ to be empty).
Let the function $\gamma(x)$ be chosen to be
\bea
\gamma(x)=A_2 x-{1\over x+C_4}
\eea
where $C_4=-{\rm inf}_M \vp+1$, and $A_2$ will be chosen later.
Note that in the range $x+C_4\geq 1$
\bea
&&
A_2x -1\leq \gamma(x)\leq A_2x, \nonumber\\
&&
A_2\leq \gamma'(x) \leq A_2+1 \nonumber\\
&&
\gamma''(x)=-2{1\over (x+C_4)^3}.
\eea
Let $\alpha(z)$ be the corresponding function as in (\ref{alpha}). If we choose $A_2\geq B$ then
by hypothesis, the function $\alpha(z)$ attains its maximum
somewhere at a point $p$ in $\bar M\setminus S$.
It suffices to show that there is a constant $C_4$
depending only on the quantities in (\ref{constants}) so
that
\bea
\alpha(p)\leq C_5.
\eea
It follows that for any $z\in \bar M\setminus S$,
\bea
\log |\nabla \vp(z)|^2-A_2 \vp(z)\leq C_5 +{1\over \vp(z)+C_4}
\leq C_6,
\eea
which implies the desired estimate.

\medskip
If $p$ is on $\pl M$, $\alpha(p)$ is immediately bounded by constants
of the form (\ref{constants}), and we are done.

\medskip
Assume then that $p$ is an interior maximum point.
Then $\Delta'\alpha(p)\leq 0$.
We can
assume that $\beta(p)\geq 1$, otherwise
\bea
\alpha(p)=\log\, \beta(p)-A_2 \vp(p)
\leq A_2 (-{\rm inf}_M\vp)
\leq C_6
\eea
and we are again done. We apply now Blocki's identity in the form (\ref{Blocki1}),
and simplify the right hand side by
dropping the terms $|\nabla'\bar\nabla \vp|^2$, ${\gamma'(\vp)\over\beta}$
on the right hand side,
\bea
0\geq\Delta'\alpha
\geq (A_2-\Lambda-F_1'-F_1''){\rm Tr}\,h^{-1}
-\gamma''(\vp)|\nabla'\vp|^2
- C_7.
\eea
Choose $A_2={\rm max}(B,\Lambda+F_1'+F_1''+1)$.
This implies that ${\rm Tr}\,h^{-1}$ is bounded above, and hence $\lambda_i^{-1}$
are all bounded above, where $\lambda_i$ are the eigenvalues of $h$.
By the Monge-Amp\`ere equation, the product of the $\lambda_i$ is bounded
by ${\rm sup}_MF$. Thus the eigenvalues $\lambda_i$ are also bounded above.
This implies that
\bea
|\nabla \vp(p)|^2\leq C_8 |\nabla' \vp(p)|^2.
\eea
It follows from the previous inequality and the explicit expression for $\gamma''$ that
\bea
{1\over (\vp+C_4)^3}|\nabla \vp|^2 \leq C_9,
\eea
that is
\bea
|\nabla \vp(p)|^2 \leq C_9 (\vp(p)+C_4)^3.
\eea
Now we may also assume that
\bea
\alpha(p)\geq 0
\eea
since otherwise there is nothing to prove. But then
\bea
A_2 \vp(p)-1\leq \gamma(\vp(p))\leq \log\,|\nabla \vp(p)|^2
\eea
and thus
\bea
\vp(p) \leq C_{10}\log\,|\nabla \vp(p)|^2+C_{11}.
\eea
Altogether we obtain
\bea
|\nabla \vp(p)|^2 \leq C_9(C_{10}\log |\nabla \vp(p)|^2+C_{12})^3.
\eea
This shows that $|\nabla \vp(p)|^2\leq C_{13}$, and since we still have
\bea
\alpha (p)
=
\log |\nabla \vp(p)|^2-\gamma(\vp(p))
\leq
\log\, C_{13}+C_{14}=C_{15}.
\eea
The theorem is proved.

\section{Application to an observation of Tsuji}
\setcounter{equation}{0}

We illustrate the application of Theorem \ref{main} to the construction
of geodesic rays from a test configuration, as in \cite{PS09}.
In the geodesic problem, we encounter a Dirichlet problem for a Monge-Amp\`ere equation
with a background $(1,1)$-form $\o_0$ which is closed, but which
may be degenerate (more precisely, strictly positive but with no positive
lower bound),
\bea
(\o_0+{i\over 2}\ddb \vp)^n= G(z)\,\o_0^n, \qquad \vp=\vp_b\ {\rm on}\ \partial M.
\eea
The degeneracy of $\o_0$ prevents the application of the standard
gradient estimates. However, in the situation of \cite{PS09},
there is a number $\kappa>0$ and an effective divisor $E$ disjoint from $\partial M$,
with $\omega_0-\kappa[E]>0$, that is, $O(E)$ admits a metric $H(z)$ so that
\bea
\omega\equiv\omega_0+{i\over 2}\ddb\log H(z)^\kappa>0.
\eea
Let $\sigma(z)$ be the canonical section of $O(E)$ which vanishes exactly on $E$.
Then a useful observation going back to Tsuji \cite{T} is that the original equation
can now be re-written as a new equation with non-degenerate background,
\bea
\label{shiftedequation}
(\o+{i\over 2}\ddb \psi)^n=F\o^n,
\eea
where $F\equiv G(\o_0^n/\o^n)$, and the new unknown $\psi$ is defined by
\bea
\psi(z)=\vp(z)-\log\,\|\sigma(z)\|^\kappa
\eea
with $\|\sigma(z)\|^2\equiv |\sigma(z)|^2H(z)$ the square of the norm of $\sigma(z)$
with respect to the metric $H(z)$. If $\vp(z)$ is bounded from below,
then $\psi(z)$ is bounded from below, and
Theorem \ref{main} applied to the shifted equation (\ref{shiftedequation}) gives upper bounds for
$|\na\psi(z)|$, and thus for $|\na\vp(z)|$ away from the divisor $E$. Precise statements
are given in Theorems 1-3 of \cite{PS09}. The main implication is that the geodesic rays
associated to general test configurations by Bergman approximations \cite{PS06, PS07, PS07a}
are $C^{1,\alpha}$ for any $0<\alpha<1$. Other related $C^{1,\alpha}$ regularity results for
geodesic rays or solutions of the degenerate Monge-Amp\`ere equation in various geometric
situations can be found in \cite{PS07a, C08, CT, BD, SZ06, SZ08}.

\section{Related estimates and other versions}
\setcounter{equation}{0}

This section is devoted to a few simple remarks and extensions of
the previous gradient estimates.

\subsection{The $C^2$ estimate}

We observed earlier that the gradient estimates of Theorem \ref{main} a)
are  natural analogues of the classic estimates of Yau
and Aubin \cite{Y78,A} for the Laplacian $\Delta\vp$
of $\vp$. This analogy
extends to the case (b) as well: assume now that $\vp(z)\to +\infty$
as $z\to S$, and that there exists a constant $C$ so that
\bea
\log(n+\Delta\vp(z))-C\vp(z)\to -\infty
\ \ {\rm as}\ \ z\to S.
\eea
Then the same a priori estimate as in (\ref{C2equation}) holds,
for $z\in \bar M\setminus S$.
This is proved by the same argument as in Yau and Aubin,
and is familiar to experts in the field. See related estimates
in e.g. \cite{EGZ, DP, TZ, ST, ST09, T}, where the preliminary $C^0$ estimates
are built on those of Kolodziej \cite{K98}.

\subsection{A gradient estimate for general Hermitian backgrounds}

The same argument applies to the case of $\omega$
just a Hermitian metric which is not necessarily K\"ahler.
This is because Blocki's identity can also be adapted to
the Hermitian case
(see e.g. eqs. (2.7)-(2.8) of \cite{PS09})
if we replace the lower bound for the bisectional curvature
by the constant $\Lambda^H$ defined by
\bea
M^{k\bar l}
R^p{}_{k\bar l}{}^{\bar m}
N_{\bar m p}
\geq -\Lambda^H ({\rm Tr}\,M)({\rm Tr} N)
\eea
for all Hermitian non-negative matrices $M$ and $N$. We can then write
\bea
(g')^{k\bar l}\pl_p u\,R^p{}_{k\bar l}{}^{\bar m}\pl_{\bar m}u
\geq -\Lambda^H \,\beta\,{\rm Tr}\, h^{-1}.
\eea
and obtain the exact same inequality as (\ref{Blocki}),
with $\Lambda$ replaced by $\Lambda^H$.
The rest of the proof proceeds as before.

We note that there has been considerable progress recently
in the study of the Monge-Amp\`ere equation on general Hermitian manifolds
\cite{GL, GLa, TW, DK}.

\bigskip
\centerline{}
\bigskip

D.H. Phong

Department of Mathematics, Columbia University, New York, NY 10027

\bigskip

Jacob Sturm

Department of Mathematics, Rutgers University, Newark, NJ 07102


\begin{thebibliography}{99}

{\small

\bibitem[A]{A} Aubin, T. ``Equations du type Monge-Amp\`ere sur les vari\'et\'es K\"ahleriennes compacts", Bull. Sc. Math. {\bf 102} (1976), 119-121 


\bibitem[BD] {BD} Berman, R. and J.P. Demailly,
``Regularity of plurisubharmonic upper envelopes
in big cohomology classes", arXiv:0905.1246

\bibitem[B1]{B08} Blocki, Z.,
``A gradient estimate in the Calabi-Yau theorem",
Math. Ann. {\bf 344} (2009), 317-327.

\bibitem[B2]{B} Blocki, Z.,
``On geodesics in the space of K\"ahler metrics",
2009 preprint.

\bibitem[C1] {C00} Chen, X.X.,
``The space of K\"ahler metrics'', J. Differential Geom.
{\bf 56} (2000) 189-234.

\bibitem[C2] {C08} Chen, X.X.,
``Space of K\"ahler metrics III: on the lower bound of the Calabi
energy and geodesic distance'', arXiv: math.DG / 0606228.

\bibitem[CT] {CT} Chen, X.X. and Y. Tang,
``Test configurations and geodesic rays'',
arXiv:0707.4149

\bibitem[DP] {DP} Demailly, J.P. and M. Pali,
``Degenerate complex Monge-Amp\`ere equations
over compact K\"ahler manifolds", arXiv: math.DG/0710.5109

\bibitem[DK] {DK} Dinew, S. and S. Kolodziej,
``Pluripotential estimates on compact Hermitian manifolds",
arXiv:0910.3937

\bibitem[EGZ] {EGZ} Eyssidieux, P., V. Guedj, and A. Zeriahi,
``Singular K\"ahler-Einstein metrics", arXiv: math/0603341

\bibitem[GL] {GL} Guan, B. and Q. Li,
``Complex Monge-Amp\`ere equations on Hermitian manifolds"
arXiv:0906.3548

\bibitem[GL1] {GLa} Guan, B. and Q. Li,
``Complex Monge-Amp\`ere equations and totally real submanifolds",
arXiv:0910.1851


\bibitem[G] {G} Guan, P.,
``A gradient estimate for complex Monge-Amp\`ere equation",
2008 preprint.


\bibitem[H] {Hanani} Hanani, A.,
``Equations du type de Monge-Amp\`ere sur les varietes hermitiennes compactes",
J. Functional Anal. {\bf 137} (1996) 49-75.


\bibitem[K]{K98} Kolodziej, S.,
``The complex Monge-Amp\`ere equation",
Acta Math. {\bf 180} (1998) 69-117


\bibitem[PS1] {PS06} Phong, D.H. and J. Sturm,
``The Monge-Amp\`ere operator and geodesics in the space of
K\"ahler potentials'', Invent. Math. {\bf 166} (2006) 125-149,
arXiv: math/0504157

\bibitem[PS2]{PS07} Phong, D.H. and
J. Sturm,
``Test configurations and geodesics in the space of K\"ahler potentials'',
J. Symplectic Geom.  {\bf 5}  (2007),  no. 2, 221-247,
arXiv: math/0606423

\bibitem[PS3]{PS07a} Phong, D.H. and J. Sturm,
``On the $C^{1,1}$ regularity of geodesics defined by test
configurations'', arXiv: math/07073956 [math.DG]

\bibitem[PS4]{PS09} Phong, D.H. and J. Sturm,
``The Dirichlet problem for degenerate complex Monge-Amp\`ere equations'',
arXiv: 0904.1898

\bibitem[PS5]{PS09a} Phong, D.H. and J. Sturm,
``Regularity of geodesic rays and Monge-Ampere equations",
arXiv:0908.0556


\bibitem[ST1] {ST} Song, J. and G. Tian,
``The K\"ahler-Ricci flow on surfaces of positive Kodaira dimension",
Inventiones Math. {\bf 170} (2007) no. 3, 609-653

\bibitem[ST2] {ST09} Song, J. and G. Tian,
``The Kahler-Ricci flow through singularities",
arXiv:0909.4898

\bibitem[SZ1] {SZ06} Song, J. and S. Zelditch,
``Bergman metrics and geodesics in the space of K\"ahler metrics on toric varieties'',
arXiv:0707.3082

\bibitem[SZ2] {SZ08} Song, J. and S. Zelditch,
``Test configurations, large deviations and geodesic rays on toric varieties'',
arXiv:0712.3599




\bibitem[T] {T} Tosatti, V.,
``Adiabatic limits of Ricci-flat Kahler metrics",
arXiv:0905.4718, to appear in J. Differential Geometry.


\bibitem[TW] {TW} Tosatti, V. and B. Weinkove,
``The complex Monge-Ampère equation on compact Hermitian manifolds",
arXiv:0910.1390


\bibitem[T1]{Tsuji} Tsuji, H.,
``Existence and degeneration of K\"ahler-Einstein metrics
on minimal algebraic varieties of general type",
Math. Ann. 281 (1988), no. 1, 123-133

\bibitem[TZ]{TZ} Tian, G. and Z. Zhang,
``On the K\"ahler-Ricci flow on projective manifolds
of general type", Chinese Ann. Math. Series B 27 (2006), no. 2, 179-192

\bibitem[Y] {Y78} Yau, S.T.,
``On the Ricci curvature of a compact K\"ahler manifold
and the complex Monge-Ampere equation I'',
Comm. Pure Appl. Math. {\bf 31} (1978) 339-411
(1978)










}

\end{thebibliography}
\end{document}